\documentclass[a4paper,12pt]{amsart}
\usepackage{color}
\usepackage{amsmath,amssymb,amscd,amsfonts,amsthm,stmaryrd}
\usepackage[mathscr]{eucal}
\usepackage{xcolor}
\usepackage{graphicx}
\usepackage{tikz-cd}
\usepackage{adjustbox}

\usepackage{nicefrac}

\usepackage{color}

\numberwithin{equation}{section}

\def\diam{\operatorname{diam}}

\def\acts{\curvearrowright}
\def\D{\partial}

\def\R{\mathbb R}
\def\Z{\mathbb Z}

\def\N{\mathbb N}

\def\ka{\kappa}

\def\eps{\epsilon}
\def\ga{\gamma}
\def\Ga{\Gamma}

\def\la{\lambda}
\def\La{\Lambda}

\def\si{\sigma}
\def\Si{\Sigma}

\def\tits{\partial_{T}}
\def\geo{\partial_{\infty}}

\def\acts{\curvearrowright}

\def\Isom{\mathop{\hbox{Isom}}}
\def\Stab{\mathop{\hbox{Stab}}}

\def\<{\langle}
\def\>{\rangle}

\theoremstyle{plain}
\newtheorem{thm}{Theorem}[section]
\newtheorem{lem}[thm]{Lemma}
\newtheorem{prop}[thm]{Proposition}
\newtheorem{cor}[thm]{Corollary}
\newtheorem{slem}[thm]{Sublemma}

\newtheorem{introthm}{Theorem}

\newtheorem{introcor}[introthm]{Corollary}

\newtheorem*{conj}{Conjecture}

\newtheoremstyle{named}{}{}{\itshape}{}{\bfseries}{.}{.5em}{\thmnote{#3} #1}
\theoremstyle{named}

\theoremstyle{definition}
\newtheorem{dfn}[thm]{Definition}

\theoremstyle{remark}

\newcommand{\bcl}{\begin{claim}}
\newcommand{\ecl}{\end{claim}}
\newcommand{\bcor}{\begin{cor}}
\newcommand{\ecor}{\end{cor}}
\newcommand{\bdfn}{\begin{dfn}}
\newcommand{\edfn}{\end{dfn}}
\newcommand{\ben}{\begin{enumerate}}
\newcommand{\bit}{\begin{itemize}}
\newcommand{\blem}{\begin{lem}}
\newcommand{\bslem}{\begin{slem}}
\newcommand{\bprop}{\begin{prop}}
\newcommand{\bthm}{\begin{thm}}
\newcommand{\een}{\end{enumerate}}
\newcommand{\eit}{\end{itemize}}
\newcommand{\elem}{\end{lem}}
\newcommand{\eslem}{\end{slem}}
\newcommand{\eprop}{\end{prop}}
\newcommand{\ethm}{\end{thm}}

\usepackage{hyperref}

\begin{document}

\title[Higher rank]{CAT(0) spaces of higher rank I}
\author{Stephan Stadler}

\newcommand{\Addresses}{{\bigskip\footnotesize
\noindent Stephan Stadler,
\par\nopagebreak\noindent\textsc{Max Planck Institute for Mathematics, Vivatsgasse 7, 53111 Bonn, Germany}
\par\nopagebreak
\noindent\textit{Email}: \texttt{stadler@mpim-bonn.mpg.de}

}}



\begin{abstract}
A CAT(0) space has rank at least $n$ if every geodesic lies in an $n$-flat.
Ballmann's Higher Rank Rigidity Conjecture predicts that a CAT(0) space of rank at least $2$ 
with a geometric group action
is {\em rigid} --
isometric to a Riemannian symmetric space, a Euclidean building, or splits as a metric product.
This paper is the first in a series motivated by Ballmann's conjecture.
Here we prove that a CAT(0) space of rank at least $n\geq 2$ is rigid if it contains a periodic $n$-flat and 
its Tits boundary has dimension $(n-1)$. This does not require a geometric group action.
The result relies essentially on the study of flats which do not bound flat half-spaces -- so-called {\em Morse flats}. 
We show that the Tits boundary $\tits F$ of a periodic Morse $n$-flat $F$ contains a {\em regular point} --
a point with a Tits-neighborhood entirely contained in $\tits F$. More precisely, we show that the set of singular
points in $\tits F$ can be covered by finitely many round spheres of positive codimension. 
\end{abstract}

\maketitle

\tableofcontents

%

\section{Introduction}
\subsection{Main results}
A {\em Hadamard manifold} is a complete simply connected Riemannian manifold of non-positive sectional curvature.
It has {\em higher rank} if every complete geodesic line is contained in  an isometrically embedded Euclidean space of dimension at least two, a so-called {\em flat}.
The structure of manifolds of non-positive curvature and higher rank was clarified in the eighties, in the series of papers \cite{BBE_structure,BBS_structure},
culminating in Ballmann's celebrated Higher Rank Rigidity Theorem \cite{B_higher,BS_higher,EH_diff}: If a higher rank Hadamard manifold admits quotients of finite volume, then 
it is either a Riemannian symmetric space or splits as a metric product. 
So in the presence of enough symmetry, higher rank implies special geometry. In particular, the theorem identifies higher rank as the exceptional case. 

The modern theory of non-positive curvature focuses on metric spaces instead of manifolds, where the curvature condition is expressed by triangle comparison. The definition is due to Alexandrov, already in the early fifties, and later gained considerable momentum from Gromov's seminal paper \cite{G_hyp}.
Nowadays,
CAT(0) spaces -- the synthetic versions of Hadamard manifolds -- play an important role in mathematics beyond geometry and topology.
They appear in many different fields such as geometric group theory,
representation theory, arithmetic and optimization.
For an account of the huge amount of literature on metric spaces with synthetic curvature bounds, see the bibliographies in \cite{AKP, Ballmann, BH}. 
The following conjecture, formulated by Ballmann and Buyalo in \cite{BaBu_periodic}, predicts a generalization of the Higher Rank Rigidity Theorem from
Hadamard manifolds to CAT(0) spaces. It is the main motivation for the present paper.

\begin{conj}[Higher Rank Rigidity]\label{conj_rr}
Let $X$ be a locally compact geodesically complete CAT(0) space. Suppose that $\Ga$ is a group of isometries of $X$
with limit set $\La(\Ga)=\geo X$. Then $\diam(\tits X)=\pi$ implies that $X$ is a Riemannian symmetric space or a Euclidean
building of rank at least 2, or that $X$ non-trivially splits as a metric product.
\end{conj}
The condition on the diameter expresses higher rank, as it is equivalent to saying that any complete geodesic bounds a flat half-plane, cf. Section~\ref{sec_cato}.
If a group $\Ga$ acts on a CAT(0) space $X$, then its {\em limit set} $\La(\Ga)$ is the set of accumulation points  
of a $\Ga$-orbit in the ideal boundary $\geo X$. The main example of group actions on CAT(0) spaces with full limit set are {\em geometric actions}, meaning
cocompact and properly discontinuous actions by isometries.

For singular spaces, Higher Rank Rigidity is wide open. It has only been established in the following short list of special cases:

\begin{itemize} 
	\item for a 2-dimensional simplicial complex with a piecewise
smooth CAT(0) metric and a geometric group action \cite{BaBr_orbi};
\item for a 3-dimensional simplicial complex with a piecewise Euclidean CAT(0) metric and a geometric group action \cite{BB_rr};
\item for a finite dimensional CAT(0) cube complex with a geometric group action 
\cite{CS_rr}\footnote{For cube complexes, the statement of Higher Rank Rigidity takes the form of a splitting result, since
neither irreducible symmetric spaces nor irreducible Euclidean buildings are cube complexes. 
They actually do not carry a piecewise cubical CAT(0) metric at all \cite{Leeb}.};
\item for a locally compact and geodesically complete CAT(0) space such that its full isometry group does not fix a point at infinity and
which supports a cocompact isometric action by an amenable locally compact group \cite{CM_bieber};
\item for a locally compact and geodesically complete CAT(0) space with 1-dimensional Tits boundary which supports a geometric group action $\Ga\acts X$
such that the induced action on the ideal boundary is not minimal, i.e. $\geo X$ contains a proper closed $\Ga$-invariant subset \cite{Ricks_1D}.
\end{itemize}

Note that the interesting case in the next to last item concerns non-discrete groups since the only locally compact geodesically complete CAT(0)
space which admits a cocompact isometric action by an amenable discrete group is the flat Euclidean space \cite[Corollary~C]{AB_amenable}.
\medskip

\begin{introthm}\label{thm_mainA}
Let $X$ be a locally compact  CAT(0) space whose Tits boundary has dimension $n-1\geq 1$. 
Suppose that every geodesic in $X$ lies in an $n$-flat. 
If $X$ contains a periodic $n$-flat, then $X$ is a Riemannian symmetric space or a Euclidean building, or  $X$ splits as a metric product.
\end{introthm}

If the isometry group of $X$ acts cocompactly, then
the condition on the dimension of the Tits boundary simply means that $X$ does not contain an $(n+1)$-flat \cite[Theorem~C]{Kleiner}. Thus, we obtain:

\begin{introcor}\label{cor_mainB}
Let $X$ be a locally compact  CAT(0) space whose 
isometry group acts cocompactly.
Suppose that every geodesic in $X$ lies in an $n$-flat but $X$
does not admit $(n+1)$-flats. 
If $X$ contains a periodic $n$-flat, then $X$ is a Riemannian symmetric space or a Euclidean building, or  $X$ splits as a metric product.
\end{introcor}

Another consequence of Theorem~\ref{thm_mainA}
is the following.

\begin{introcor}\label{cor_mainA}
Let $X$ be a locally compact  CAT(0) space whose Tits boundary has dimension $n-1\geq 1$.
Suppose that every geodesic in $X$ lies in an $n$-flat. 
If $X$ contains a periodic $n$-flat, then $X$ has a unique decomposition into irreducible factors
\[X=\R^j\times X_1\times\ldots \times X_k\times Y_1\times\ldots \times Y_l\times M_1\times\ldots\times M_m\]
where $X_a$ are rank 1 spaces, $Y_b$ are irreducible Euclidean buildings of higher rank and $M_c$ are
irreducible Riemannian symmetric spaces of higher rank.
\end{introcor}

Note that in the smooth case, the higher rank property can be expressed by asking that every complete geodesic admits $n$ linearly independent parallel Jacobi fields for some $n\geq 2$.
In the presence of a geometric group action these Jacobi fields can then be ``integrated'' to  $n$-flats \cite[Section~IV.4]{ballmannbook}.
Hence, in the context of metric spaces, the assumption on the existence of $n$-flats in Theorem~\ref{thm_mainA} is natural\footnote{
See \cite{OR_int} for a definition of rank for CAT(0) spaces whose isometry groups satisfy the duality condition.}.
Recall that a flat  is {\em periodic}, if its stabilizer in the isometry group of the surrounding space contains a subgroup acting geometrically on it. 
By \cite[Flat Torus Theorem~7.1]{BH}, the existence of  a periodic $n$-flat in $X$ is equivalent to the presence of a subgroup $\Z^n<\Isom(X)$.
We point out that Theorem~\ref{thm_mainA} does not assume any symmetries besides the existence of a single periodic $n$-flat.
Also, if $X$ is a rank $n$ symmetric space or  Euclidean building, and $\Ga$ acts cocompactly on $X$, then $X$
contains a periodic $n$-flat \cite[Theorem~8.9]{BaBr_orbi}, \cite[Lemma~8.3]{Mos_rigidity}, \cite[Theorem~2.8]{PR_cartan}.
We emphasize that some symmetry is required in order for the conclusion of Higher Rank Rigidity to hold true, since a geodesically complete CAT(0) space whose Tits boundary $\tits X$ has
diameter $\pi$ 
does not have to be a product or a Euclidean building \cite{BB_rr}.
 \medskip

We now discribe the main novelty of the present paper.
Recall that a point in a metric space is called {\em regular}, if it has a neighborhood
homeomorphic to an open set in a Euclidean space.
The crucial ingredient behind Theorem~\ref{thm_mainA} is a contribution 
 to the following general conjecture which has been completely open.

\begin{conj}[Lytchak's Regular Point Conjecture]
Let $X$ be a locally compact CAT(0) space with a geometric group action. 
Then $\tits X$ contains a regular point.
\end{conj}

Gromov first raised the question which CAT(1) spaces do contain regular points \cite[p.90]{G_asym}.
A general CAT(1) space does not have to contain any regular points.
On the other hand, locally compact and geodesically complete CAT(1) spaces always do, 
as shown by Lytchak-Nagano in \cite{LN_gcba},
following ideas of Burago-Gromov-Perelman \cite{BGP}. However, the 
Tits boundary of a locally compact CAT(0) space with a geometric group action
is typically not locally compact.

Recall that a flat in a CAT(0) space is called {\em Morse} if it does not bound a flat half-space.
We prove:

\begin{introthm}{[Theorem~\ref{thm_finitesing}]}\label{thm_mainF}
Let $X$ be a locally compact CAT(0) space.  Suppose that $X$ contains a periodic Morse flat $F$. Then 
the complement in $\tits F$ of the set of regular points can be covered by a finite set of round spheres of positive codimension.
\end{introthm}

In particular, the singular set in the Tits boundary of a periodic Morse $(n+1)$-flat is given by a union of finitely many hyperspheres and
a closed set of finite $(n-2)$-dimensional Hausdorff measure. In this sense the Tits boundary of a periodic Morse flat resembles a Coxeter complex.
This will be used in \cite{St_rrII} in order to control the geometry of certain submetries of the Tits boundary.
Let us emphasize again that Theorem~\ref{thm_mainF} does not assume any symmetry besides the periodic Morse flat,  no geometric group action is required.

\medskip

The present article is the first in a series motivated by the Higher Rank Rigidity Conjecture.
The sequels all rely on the crucial new idea on how to find regular points at infinity presented here.  
In \cite{St_rrII}
we prove the following rigidity result which generalizes the Higher Rank Rigidity Theorem for Hadamard manifolds.

\bthm[{\cite[Main Theorem]{St_rrII}}]\label{thm_rrII}
Let $X$ be a locally compact CAT(0) space  with a geometric group action $\Ga\acts X$.  
Suppose that there exists $n\geq 2$ such that every geodesic in $X$ lies in an $n$-flat. 
If $X$ contains a periodic Morse $n$-flat, then $X$ is a Riemannian symmetric space or a Euclidean building or  $X$
non-trivially splits as a metric product.
\ethm

In \cite{St_rrIII} we confirm Higher Rank Rigidity as well as the Closing Lemma for spaces without 3-flats.
More precisely, we obtain

\bthm[{\cite[Theorems~A and~B]{St_rrIII}}]\label{thm_rrIII}
Let $X$ be a locally compact geodesically complete  CAT(0) space without 3-flats. 
Suppose $X$ admits a geometric group action $\Ga\acts X$. If $X$ contains a complete geodesic which does not bound a flat half-plane, then it also contains a
 $\Ga$-periodic geodesic which does not bound a flat half-plane. If every complete geodesic in $X$ bounds a flat half-plane,
then $X$ is a Riemannian symmetric space, a Euclidean building or non-trivially splits as a metric product.
\ethm

%
%
%
%

\subsection{Acknowledgments}
 I would like to use this opportunity to thank Bernhard Leeb and Alexander Lytchak for 
pointing out several mistakes in a late version of this article.
Their insightful comments led to substantial improvements. I also want to thank Bruce Kleiner for several inspiring discussions.
I was supported by DFG grant SPP 2026.


\section{Preliminaries}\label{sec_pre}

\subsection{Metric spaces}

We assume that the reader is familiar with the geometry of metric spaces with upper curvature bound,
as references we mention \cite{AKP, Ballmann, BH, KleinerLeeb}.
For the required preliminaries on buildings we refer the reader to \cite{KleinerLeeb,Leeb}.

Here we only want to agree on notation and collect some basic facts that will be used later.

Euclidean $n$-space with its flat metric will be denoted by $\R^n$. The unit sphere $S^{n-1}\subset\R^n$ equipped with the induced metric will be referred to as a
{\em round sphere}. Its intersection with a half-space $\R^{n-1}\times[0,\infty)$ is called a {\em round hemisphere}.
We denote  the distance between two points $x$ and $y$ in a metric space $X$ by $|x,y|$.
For $x\in X$ and $r>0$, we denote by $B_r(x)$ and $\bar B_r(x)$ the open and closed $r$-ball around $x$, respectively.
Similarly, $N_r(A)$ and $\bar N_r(A)$ denote the open and closed $r$-neighborhood of a subset $A\subset X$, respectively.
A \emph{geodesic}
is an isometric embedding of an interval. It is called a {\em geodesic segment}, if it is compact.
The {\em endpoints} of a geodesic segment $c$ are denoted by $\D c$.

A {\em branch point} is a point through which a certain geodesic can be extended in more than one way.
More precisely, a point $b$ is a branch point if there are geodesic segments $c$, $c^+$ and $c^-$
such that $c\cup c^\pm$ is a geodesic and $c^+\cap c^-=\{b\}$.
In a \emph{geodesic metric space} 
any pair of   points
is connected by a geodesic.
Such a space is \emph{geodesically complete} if every geodesic segment is contained in a complete local geodesic.

\subsection{Spaces with an upper curvature bound}

For any CAT($\kappa$) space  $X$,
the angle between each pair of geodesics starting at the same point
is well defined. 
 Let $x,y,z$ be three points at pairwise distance less than $\frac{\pi}{\sqrt{\kappa}}$ in a CAT($\kappa$) space $X$.
Whenever $x\neq y$, the geodesic between $x$ and $y$ is unique and will be denoted
by $xy$.   For $y,z \neq x$, the angle at $x$ between $xy$ and $xz$
will be denoted by $\angle_x(y,z)$.
The \emph{space of directions} or {\em link}
 at a point  $x\in X$ will be denoted by $(\Si_x X,\angle)$,
its elements are called {\em directions (at $x$)}.

The {\em dimension $\dim (X)$} of a CAT($\ka$) space $X$ refers to Kleiner's 
geometric dimension \cite{Kleiner}.
It vanishes precisely when the space is discrete.
In general, it is defined inductively:
\[\dim (X)= \sup _{x\in X} \{ \dim (\Sigma _xX) +1 \}.\]

\subsection{CAT(1) spaces}
  
For two CAT(1) spaces $Z_1$ and $Z_2$ we denote by $Z_1\circ Z_2$ their {\em spherical join}.
It is a CAT(1) space of diameter $\pi$. A pair of points $\xi^-,\xi^+$ in a CAT(1) space $Z$ are called {\em antipodes},
if their distance is at least $\pi$.

A subset  $A\subset Z$ is called {\em spherical}, if it 
embeds isometrically into a round sphere. 
A {\em (spherical) $n$-lune of angle $\theta\in[0,\pi]$} in a CAT(1) space is a closed convex subset $\la$ isometric to 
$S^{n-1}\circ[0,\theta]$.
Bigons in CAT(1) spaces lead to spherical lunes \cite[Lemma~2.5]{BB_diam}.
We will use the following more general result which follows immediately from \cite[Lemma~2.5]{BB_diam} and \cite[Lemma~4.1]{Ly_rigidity}.

\blem[Lune Lemma]\label{lem_lune}
Let $\xi$ and $\hat\xi$ be antipodes at distance $\pi$ in a CAT(1) space $Z$.
If $\ga^\pm$ are geodesics from $\xi$ to $\hat\xi$, then
$\ga^-\cup\ga^+$ is a spherical subset in $Z$. 
More precisely, if $\angle_\xi(\ga^-,\ga^+)\geq \pi$, then $\ga^-\cup\ga^+$ is a round $1$-sphere in $Z$;
and if  $\angle_\xi(\ga^-,\ga^+)=\theta< \pi$, then $\ga^-\cup\ga^+$ bounds a spherical lune of angle $\theta$.
Similarly, let $\tau^\pm\subset Z$ be round $n$-hemispheres with $\D\tau^-=\D\tau^+$ and let $\zeta^\pm\in Z$
denote their centers. If $|\zeta^-,\zeta^+|\geq\pi$, then $\tau^-\cup\tau^+$ is a round $n$-sphere in $Z$;
and if $|\zeta^-,\zeta^+|=\theta<\pi$, then $\tau^-\cup\tau^+$ bounds a spherical $(n+1)$-lune of angle $\theta$.
\elem


\bdfn\label{def_reg_point}
Let $Z$ be a CAT(1) space. We call a point $\xi\in Z$ {\em $k$-regular}, if
it has a neighborhood homeomorphic to an open set in $\R^k$. The point $\xi$ is called {\em $k$-spherical},
if there exists  a radius $r>0$ such that $B_r(\xi)\subset Z$ is isometric to an $r$-ball in $S^{k}$.
In both definitions, if it's clear from the context, we neglect the quantification and speak of {\em regular} respectively {\em spherical} points.
\edfn
%

\blem[{\cite[Lemma~2.1]{BL_building}}]\label{lem_antintop}
	Let $Z$ be a CAT(1) space of dimension $n$. If $\si\subset Z$ is a round $n$-sphere, then every point $\xi\in Z$ has an antipode in $\si$.
\elem


This implies the following well-known criterion for geodesic completeness.

\blem\label{lem_gc}
	Let $Z$ be a CAT(1) space of dimension $n$. If every point $\xi\in Z$ is contained in a round $n$-sphere in $Z$,
	then $Z$ is geodesically complete.
\elem

\proof
If $\xi$ is a point in $Z$ contained in a round $n$-sphere $\si\subset Z$, 
then any local geodesic ending in the point $\xi$ can be extended beyond  $\xi$
by a geodesic in $\si$, since any direction in $\Si_\xi Z$ has an antipode in $\Si_\xi\si\cong S^{n-1}$ by Lemma~\ref{lem_antintop}.
\qed
\medskip

A simple but useful property in CAT(1) spaces with a geodesically complete link is the following.

\blem\label{lem_specialanti}
Let $Z$ be a CAT(1) space and $z\in Z$ a point where the link $\Si_zZ$ is geodesically complete. Then for every pair of directions $v, w\in\Si_z Z$
there exists an antipode $\hat w\in\Si_z Z$ of $w$, such that 
\[\angle_z(\hat w,v)=\pi-\angle_z(w,v).\]
\elem

\proof
If $\angle_z(v,w)=\pi$, then we can take $\hat w=v$. Otherwise, we extend the geodesic from $w$ to $v$ up to distance $\pi$.
The endpoint $\hat w$ is then an antipode of $w$ with the required property. 
\qed
\medskip

The proof of Theorem~\ref{thm_mainA} relies on the following rigidity result of Lytchak.

\bthm[{\cite[Main Theorem]{Ly_rigidity}}]\label{thm_lytchak}
Let $Z$ be a finite dimensional geodesically complete
CAT(1) space. If $Z$ has a proper closed subset $A$ containing with each
point all of its antipodes, then $Z$ is
a spherical join or a spherical building.
\ethm

\subsection{CAT(0) spaces}\label{sec_cato}

The {\em ideal boundary} of a CAT(0) space $X$, equipped with the cone topology, is denoted by $\geo X$. 
If $X$ is locally compact, then $\geo X$ is compact.  
The {\em Tits boundary} of $X$ is denoted by $\tits X$, it is the ideal boundary equipped with 
the Tits metric $|\cdot,\cdot|_T$. Recall that the Tits metric is the intrinsic metric associated to the {\em Tits angle}. 

If $C$ is a closed convex subset, then it is CAT(0) with respect  to the induced metric. In this case, $C$ admits a 1-Lipschitz retraction $\pi_C:X\to C$.
If $C_1$ and $C_2$ are closed convex subsets, then the distance function $d(\cdot,C_1)|_{C_2}$ is convex; and constant if
and only if $\pi_{C_1}$ restricts to an isometric embedding on $C_2$. We call $C_1$ and $C_2$ {\em parallel}, $C_1\| C_2$, if
and only if  $d(\cdot,C_1)|_{C_2}$ and $d(\cdot,C_2)|_{C_1}$ are constant.
Let $Y\subset X$ be a geodesically complete closed convex subset. Then we define the {\em parallel set} $P(Y)$
as the union of all closed convex subsets parallel to $Y$. The parallel set is closed, convex and splits canonically as a metric product
\[P(Y)\cong Y\times CS(Y)\]
where the {\em cross section} $CS(Y)$ is a closed convex subset.

If $X_1$ and $X_2$ are CAT(0) spaces, then their metric product $X_1\times X_2$ is again a CAT(0) space.
We have $\tits (X_1\times X_2)=\tits X_1\circ \tits X_2$ and $\Si_{(x_1,x_2)}(X_1\times X_2)=\Si_{x_1} X_1\circ \Si_{x_2} X_2$.  
If $X$ is a geodesically complete CAT(0) space, then any join decomposition of $\tits X$ is induced by a metric product decomposition of $X$ \cite[Proposition~2.3.7]{KleinerLeeb}.

The Tits boundary of a closed convex subset $Y\subset X$ embeds canonically $\tits Y\subset\tits X$.
If two closed convex subsets $Y_1$ and $Y_2$ intersect in $X$, then 
\[\tits(Y_1\cap Y_2)=\tits Y_1\cap \tits Y_2.\]


For points $\xi\in\geo X$ and $x\in X$ we denote the {\em Busemann function centered at $\xi$ based at $x$} by $b_{\xi,x}$.
If $\rho:[0,\infty)\to X$ denotes the geodesic ray asymptotic to $\xi$ and with $\rho(0)=x$, then
\[b_{\xi,x}(y)=\lim\limits_{t\to\infty}(|y,\rho(t)|-t).\]
It is a 1-Lipschitz convex function whose negative gradient at a point $y\in X$ is given by $\log_y(\xi)$. 
We denote the {\em horoball} centered at a point $\xi\in\geo X$ and based at the point $x\in X$ by 
\[HB(\xi,x):=b^{-1}_{\xi,x}((-\infty,0]).\]
It is a closed convex subset with 
\[\tits HB(\xi,x)=\bar B_{\frac{\pi}{2}}(\xi)\subset\tits X.\]

A {\em $n$-flat} $F$ in a CAT(0) space $X$ is a closed convex subset isometric to $\R^n$.
In particular, $\tits F\subset\tits X$ is a round $(n-1)$-sphere.
On the other hand, if $X$ is locally compact and $\si\subset \tits X$ is a round  $(n-1)$-sphere, then either there exists an
$n$-flat $F\subset X$ with $\tits F=\si$, or there exists a round $n$-hemisphere $\tau^+\subset\tits X$
with $\si=\D\tau^+$ \cite[Proposition~2.1]{Leeb}. Consequently, if $\tits X$ is $(n-1)$-dimensional, then any round $(n-1)$-sphere in $\tits X$
is the Tits boundary of some $n$-flat in $X$. Moreover, a round $n$-hemisphere $\tau^+\subset\tits X$ bounds a flat $(n+1)$-half-space in $X$ if and only if
its boundary $\D\tau^+$ bounds an $n$-flat in $X$. A flat is called {\em Morse}, 
if it does not bound a flat half-space\footnote{For a motivation of the term and further context see \cite{HKS_I,HKS_II}.}.

The following basic fact will be used repeatedly throughout the paper.

\blem[{\cite[Sublemma~2.3]{Leeb}}]\label{lem_flats_in_convex}
Let $F$ be a flat in a CAT(0) space and let $C\subset X$ be a closed convex subset so that $\geo F\subset\geo C$. 
Then $C$ contains a flat $F'$ parallel to $F$.
\elem


A {\em flat ($n$-dimensional) half-space} $H\subset X$ is a closed convex subset isometric to a Euclidean half-space $\R_+^n$.
Its boundary $\D H\subset X$ is an $(n-1)$-flat and its Tits boundary is a round $(n-1)$-hemisphere $\tits H\subset\tits X$.  Flat half-spaces will play a central role in our arguments later, and we agree to denote them by $H$
and their boundaries by $\D H=h$.

We define the {\em parallel set} of a round sphere $\si\subset \tits X$ as 
\[P(\si)=P(F)\]
where $F$ is a flat in $X$ with $\tits F=\si$, if such a flat exists. 

The proof of Theorem~\ref{thm_mainA} relies on the following rigidity result of Leeb.

\bthm[{\cite[Main Theorem]{Leeb}}]\label{thm_leeb}
Let $X$ be a locally compact, geodesically complete CAT(0) space. If $\tits X$ is a connected thick irreducible spherical building,
then $X$ is either a symmetric space or a Euclidean building.
\ethm

\section{Recognizing a building at infinity}\label{sec_buildings}

The aim of this section is to prove the following result. 

\bthm\label{thm_rec}
Let $X$ be a locally compact CAT(0) space where every geodesic lies in an $n$-flat. Suppose that $\dim(\tits X)=n-1$.
If $\tits X$ contains an open
relatively compact subset $U$, then $\tits X$ is a spherical building.
\ethm

The result is similar to \cite[Theorem~1.6]{BL_building} and builds on Theorem~\ref{thm_lytchak} in the same way.
Note that if we would know that any pair of antipodes in $\tits X$ lies in a round $(n-1)$-sphere, 
then we could directly apply \cite[Theorem~1.6]{BL_building}.

Using the structure theory of locally compact and geodesically complete spaces with upper curvature bounds \cite{LN_gcba},
one can deduce Theorem~\ref{thm_rec} form the following similar result which will suffice for our proof of Theorem~\ref{thm_mainA}.

\bthm\label{thm_rec_sph}
Let $X$ be a locally compact CAT(0) space where every geodesic lies in an $n$-flat. Suppose that $\dim(\tits X)=n-1$.
If $\tits X$ contains an $(n-1)$-spherical point, then $\tits X$ is a spherical building.
\ethm

We need some preparation before we can provide the proof of Theorem~\ref{thm_rec}.

%

Recall that a pair of antipodes in the Tits boundary of a CAT(0) space $X$ does not have to bound a complete geodesic in $X$.
Nevertheless, for locally compact geodesically complete spaces we have the following.

\blem\label{lem_antipodes}
Let $X$ be a locally compact, geodesically complete CAT(0) space. Let $\xi$ and $\hat\xi$ be antipodes in $\tits X$.
Then there exists a sequence of complete geodesics $c_k$ in $X$ with $\geo c_k=\{\hat\xi,\xi_k\}$ and $\xi_k\to\xi$ in $\geo X$.
\elem

\proof
If $|\xi,\hat\xi|>\pi$, then there is a complete geodesic $c$ in $X$ with $\geo c=\{\xi,\hat\xi\}$ \cite[Proposition~2.1]{Leeb}.
Hence, we may assume $|\xi,\hat\xi|=\pi$. 
Fix a base point $o\in X$ and let $\rho:[0,\infty)\to X$ be a parametrization of the geodesic ray $o\xi$.
Since $X$ is locally compact and geodesically complete,  all links in $X$ are geodesically complete \cite[Corollary~5.8, Corollary~5.9]{LN_gcba}.
By Lemma~\ref{lem_specialanti}, we can choose for each $k\in\N$ a complete geodesic $c_k$ extending the geodesic ray $\rho(k)\hat\xi$ and such that 
\[\angle_{\rho(k)}(\hat\xi,\xi)=\pi-\angle_{\rho(k)}(\xi,\xi_k)\]
where $\xi_k:=c_k(+\infty)$. 
Now for every $p\in o\xi$ we have
\[\angle_p(\xi_k,\xi)\leq\angle_{\rho(k)}(\xi_k,\xi)=\pi-\angle_{\rho(k)}(\xi,\hat\xi)\to 0.\]
Since $\angle_{\rho(k)}(\xi,\hat\xi)\to|\xi,\hat\xi|=\pi$. As $p\in o\xi$ was arbitrary, this shows $\xi_k\to\xi$ in $\geo X$.
Indeed, by local compactness, we may assume $\xi_k\to\xi'$ after passing to a subsequence. In particular, $\angle_x(\xi_k,\xi')\to 0$ for all $x\in X$. 
Now if $\xi\neq\xi'$, then there exists $q\in o\xi$
such that $\angle_q(\xi,\xi')>0$. But 
\[\angle_q(\xi,\xi')\leq \angle_q(\xi,\xi_k)+\angle_q(\xi_k,\xi')\to 0.\]
\qed
\medskip

\blem\label{lem_anti_spheres}
Let $X$ be a locally compact CAT(0) space. Suppose that every geodesic in $X$ lies in an $n$-flat.
If $\xi\in\tits X$ is an $(n-1)$-spherical point, then for every antipode $\hat\xi$ of $\xi$ there exists a round $(n-1)$-sphere $\si\subset\tits X$
containing $\xi$ and $\hat\xi$.
\elem

\proof
Let
$\xi\in\tits X$ be an $(n-1)$-spherical point, cf. Definition~\ref{def_reg_point}.
Then there exists $s>0$ such that $B_s(\xi)$ is contained in any $(n-1)$-sphere which contains $\xi$.
If $\hat\xi\in\tits X$ is an antipode of $\xi$, then by Lemma~\ref{lem_antipodes}, there exists a sequence of complete geodesics $(c_k)$
in $X$ with $\geo c_k=\{\xi,\hat\xi_k\}$ and $\hat\xi_k\to\hat\xi$. By assumption, there exists $n$-flats $F_k$ with $c_k\subset F_k$. 
Hence, if $\eta\in\tits X$ is a point with $|\eta,\xi|=s$, then $|\eta,\hat\xi_k|=\pi-s$ for all $k\in\N$. Thus,
\[\pi=|\xi,\hat\xi|\leq|\xi,\eta|+|\eta,\hat\xi|\leq s+\liminf\limits_{k\to\infty}|\eta,\hat\xi_k|=\pi,\]
by lower semi-continuity of the Tits metric. We obtain $|\eta,\hat\xi|=\pi-s$ for every point $\eta$ at distance $s$
from $\xi$.
In particular, through every such point runs a geodesic from $\xi$ to $\hat\xi$.
Denote by $\si\subset\tits X$ the union of all these geodesics. It follows from Lemma~\ref{lem_lune}
that $\si$ is a round $(n-1)$-sphere.
\qed

\bprop\label{prop_reg_anti}
Let $X$ be a locally compact CAT(0) space where every geodesic lies in an $n$-flat.
Suppose $\dim(\tits X)=n-1$.
Then the set $O\subset \tits X$ of $(n-1)$-spherical points is closed under taking antipodes. 
\eprop

\proof
Let $\xi$ be a $(n-1)$-spherical point in $\tits X$ and let $\hat \xi$ be an antipode.
By Lemma~\ref{lem_anti_spheres}, there exists a round $(n-1)$-sphere $\si\subset \tits X$ which contains $\xi$
and $\hat\xi$. Since $\xi$ is spherical, there exists $s>0$ such that $B_s(\xi)\subset\si$.
We claim that $B_s(\hat\xi)\subset\si$ holds as well.
Let $\hat\eta$ be a point in $B_s(\hat\xi)$. Since $\dim(\tits X)=n-1$,
 we can extend the geodesic $\hat\eta\hat\xi$, as in the proof of Lemma~\ref{lem_gc},
up to an antipode $\eta$ of $\hat\eta$ in $\si$. Then $|\eta,\xi|=|\hat\eta,\hat\xi|<s$ and therefore $\eta$ is spherical.
By Lemma~\ref{lem_anti_spheres}, we find a round $(n-1)$-sphere $\si'\subset \tits X$ which contains $\eta$
and $\hat\eta$. Because $\eta\in B_s(\xi)$, we must have $B_s(\xi)\subset\si\cap\si'$. By construction,
$\hat \xi$ also lies in $\si\cap\si'$. Now convexity of $\si\cap\si'$ implies $\si=\si'$ and therefore $\hat\eta\in\si$ as required.
\qed

\blem\label{lem_splitflat}
Let $X$ be a locally compact CAT(0) space with $\dim(\tits X)=n-1$ and where every geodesic lies in an $n$-flat.
Suppose $X$ splits as a non-trivial metric product $X\cong X_1\times X_2$.
Then we have $n=n_1+n_2$ where $\dim(\tits X_j)=n_j-1$, $j=1,2$, 
and every complete geodesic in $X_j$ is contained in an $n_j$-flat.
\elem

\proof
It is clear that the rank is additive. So
let $c_j$ be a complete geodesic in $X_j$ and let $\tilde c=(c_1,c_2)$ be the diagonal geodesic in
$c_1\times c_2\subset X_1\times X_2$. If $\tilde F$ is an $n$-flat in $X$ which contains $\tilde c$,
then by \cite[Lemma~2.3.8]{KleinerLeeb} there are flats $F_j\subset X_j$ such that $\tilde F\subset F_1\times F_2$.
In particular, $\dim(F_j)=n_j$ and $c_j\subset F_j$. 
\qed
\medskip

\proof[Proof of Theorem~\ref{thm_rec_sph}]
By assumption, $\dim(\tits X)=n-1$ and every point $\xi$ in $\tits X$ lies in a round $(n-1)$-sphere. 
Hence, Lemma~\ref{lem_gc} implies that $\tits X$ is geodesically complete. 

We proceed by induction on the dimension of $\tits X$. If the $\dim(\tits X)$ is 0, then $\tits X$ is discrete and there is nothing to show. 
So let's assume $\dim(\tits X)\geq 1$.
Recall that two antipodes $\xi,\hat\xi$ in $\tits X$ at distance $>\pi$ bound a complete geodesic in $X$ \cite[Proposition~2.1]{Leeb}.
Since every complete geodesic in $X$ lies in an $n$-flat, we conclude that the diameter of $\tits X$ is $\pi$.

Suppose $X$ splits as a non-trivial metric product $X\cong X_1\times X_2$.
Then by Lemma~\ref{lem_splitflat},  the factors  are subject to the induction
hypothesis and we may assume that
$X$ is irreducible.

By assumption, the set of $(n-1)$-spherical points $O\subset\tits X$
is non-empty. If $O=\tits X$, then $X$ is isometric to $\R^n$. So let us assume that $O$ is a proper
subset of $\tits X$. By Proposition~\ref{prop_reg_anti}, the set $\tits X\setminus O$ is
a non-empty, closed proper subset which contains with every point also all its Tits antipodes.
Since $\tits X$ is irreducible, Theorem~\ref{thm_lytchak} implies that it is a spherical building.
\qed

\proof[Proof of Theorem~\ref{thm_rec}]
By Theorem~\ref{thm_rec_sph}, it is enough to find an $(n-1)$-spherical point in $\tits X$.
If $U$ is an open relatively compact subset of $\tits X$, then $U$ contains an open subset $U'$
homeomorphic to an $n$-manifold \cite[Theorem~1.2]{LN_gcba}. Hence, for a point $\xi\in U'$
there exists an $\eps>0$ such that $B_\eps(\xi)\subset U'$. If $\si\subset \tits X$
denotes a round $(n-1)$-sphere with $\xi\in\si$, then $\si\cap B_\eps(\xi)=B_\eps(\xi)$.
Hence $\xi$ is $(n-1)$-spherical.
\qed

\section{Rank $n$ spaces with periodic $n$-flats}\label{sec_rank_n}

\subsection{From branch points to orthogonal half-planes}

\bdfn\label{def_ofh}
Let $F\subset X$ be a flat in a CAT(0) space. Then we call a flat half-space $H\subset X$ a
{\em flat half-space orthogonal to $F$}, or, if the respective flat is clear from the context, an
{\em orthogonal flat half-space} if 
\[H\cap F=\D H\ \text{ and }\ \angle(H,F)\geq\frac{\pi}{2}\]
where the angle is measured in the cross section of the parallel set $P(\D H)$.
In case $H$ is 2-dimensional, we adjust to 
{\em flat half-plane orthogonal to $F$} and
{\em orthogonal flat half-plane}, respectively.
\edfn

\bdfn
Let $X$ be a locally compact CAT(0) space and $F\subset X$ a flat.
Let $\Isom(X)$ denote the isometry group of $X$.
We say that $F$ is {\em periodic}, if its {\em stabilizer} $\Stab(F):=\{\ga\in \Isom(X)|\ \ga(F)=F\}$ contains a subgroup acting geometrically on $F$.
\edfn

\blem\label{lem_bieber}
Let $X$ be a locally compact CAT(0) space and let $F\subset X$ be a periodic $n$-flat with $n\geq 2$.
Suppose that $F=F^-\times F^+$ is a metric product decomposition into flats $F^\pm$ of dimension $n^\pm$.
Then $F^\pm$ is a periodic $n^\pm$-flat in the cross section $CS(F^\mp)$.
\elem

\proof
By Bieberbach's theorem, the stabilizer of $F$ contains a free abelian subgroup $A$ of rank $n$
which acts cocompactly on $F$ by Euclidean translations. Since $A$ preserves $\tits F$ pointwise, it also preserves $P(F^\mp)$
and the splitting $P(F^\mp)\cong F^\mp\times CS(F^\mp)$.
Hence its projection to $CS(F^\mp)$ contains a free abelian subgroup of rank $n^\pm=n-n^\mp$ which  acts 
geometrically by Euclidean translations on $F^\pm\subset CS(F^\mp)$. 
\qed
\medskip

The following result implies that a branch point in the ideal boundary of a periodic flat leads to an orthogonal flat half-plane.

\blem[{\cite[Lemma~2.3.1]{HK}}]\label{lem_perpH}
Let $X$ be a locally compact CAT(0) space.
Let $F\subset X$ be a periodic flat. Let $\xi\in\tits F$ be a point such that there exists a sequence $(\xi_k)$ in $\tits X\setminus\tits F$
which converges to $\xi$ with respect to the Tits metric. 
Then there exists a flat half-plane $H\subset X$ asymptotic to $\xi$ and orthogonal to $F$ (cf. Definition~\ref{def_ofh}). 
\elem

\proof
Since $F$ is periodic, it is enough to find for every $R>0$ a flat strip $S\cong\R\times[0,R]$ which is asymptotic to $\xi$ and
orthogonal to $F$. We fix a base point $o\in F$. For every $k\in\N$ we choose $x_k\in o\xi_k$ such that $|x_k,F|=R$.
Denote by $\bar x_k\in F$ the nearest point to $x_k$. Choose isometries $\ga_k\in\Stab(F)$ such that all points $\ga_k(\bar x_k)$
lies in a fixed compact set. After passing to a subsequence, we obtain a limit geodesic $c_\infty$ which is parallel to $F$
and spans the required flats strip. More precisely, we have the convergence $\ga_k(o\xi_k,x_k)\to(c_\infty,x_\infty)$ with respect to pointed Hausdorff topology.   
The complete geodesic $c_\infty$ is asymptotic to $\xi$ since $|\ga_k\xi_k,\xi|=|\xi_k,\xi|\to 0$. Because $x_k o\subset \bar N_R(F)$, also $c_\infty$
is contained in the $R$-neighborhood of $F$ and has even distance $R$ from $F$ since $|x_\infty,F|=R$.
Denote by $\bar c_\infty\subset F$ the complete geodesic which contains the point $x_\infty$ and is asymptotic to $\xi$.
Then $c_\infty$ and $\bar c_\infty$ span the required flat strip by \cite[Lemma~2.3.5]{KleinerLeeb}.  
\qed

\bcor\label{cor_perpH}
Let $X$ be a locally compact CAT(0) space.
Let $\hat F\subset X$ be a periodic $(n+1)$-flat with $\tits \hat F=\hat\si$.
Denote by $\la\subset\tits X$  a spherical $(k+1)$-lune with $k\leq n$
bounded by round $k$-hemispheres $\tau^\pm$. Suppose that $\la\cap\hat\si=\tau^-$.
Then there exists a round $(k+1)$-hemisphere $\hat\tau\subset\tits X$
such that $\D\hat\tau\subset\hat\si$ and $\tau^-\subset\D\hat\tau$.
Moreover, $\hat\tau$ is the ideal boundary of a flat half-space $\hat H$
orthogonal to $\hat F$.
\ecor

\proof
Set $\si=\D\tau^\pm$ and write $\hat\si=\si\circ\si^\perp$.
Let $\hat F=F\times F^\perp$ be a corresponding product splitting.
We have $\la=\si\circ[0,\theta]$. By Lemma~\ref{lem_bieber}, $F^\perp$
is a periodic flat in the cross section $CS(F)$. Denote by $\zeta^-$
the center of $\tau^-$. Then the cross section of $\la$ provides a path in $\tits CS(F)$ which intersects
$\tits F^\perp$ only in  $\zeta^-$. Thus, by Lemma~\ref{lem_perpH}, there exists an orthogonal flat half-plane $H\subset CS(F)$
for $F^\perp$ which is asymptotic to $\zeta^-$. Hence, $\hat H:=F\times H\subset P(F)$ and $\hat\tau:=\tits \hat H$
define an orthogonal flat half-space and round $(k+1)$-hemisphere, respectively, as required.
\qed

\subsection{Orthogonal half-spaces from a sequence of lower dimensional ones}

In this section we prove the following key result which allows us to understand the singular set  
in the Tits boundary of a periodic Morse flat. The responsible geometric property behind the proof 
is flatness of intersecting parallel sets, similar as in \cite[Theorem~2]{St_obst}.

\bprop\label{prop_halfspace}
Let $X$ be a locally compact  CAT(0) space.
Let $\hat F$ be a periodic $\hat n$-flat in $X$ with $\tits\hat F=\hat \si$.
Suppose that there is a sequence of  flat half-spaces $(H_k)$ of dimension $d\leq n$,
all orthogonal to $\hat F$, such that their boundary flats $h_k:=\D H_k$ are  pairwise non-parallel. 
Then there exists a $(d+1)$-dimensional flat half-space $\hat H$ orthogonal to $\hat F$.
Moreover, if the sequence of round spheres $\si_k=\geo h_k$ converges to a round sphere $\si\subset\hat\si$ and for every $k\in\N$
the convex hull of $\si_k$ and $\si$ is a round sphere of dimension $\hat d-1$, then we can arrange $\hat H$ to have 
dimension $\hat d+1$ and such that $\si\subset\geo\hat H$.
\eprop

\proof
Since $\hat F$ is periodic, we may assume that we have a convergent sequence of pairwise distinct pointed flat half-spaces $(H_k,p_k)\to (H_\infty,p_\infty)$, 
where $p_k\in F_k$.
The limit $H_\infty$ with boundary flat $h_\infty:=\D H_\infty$ is then still orthogonal to $\hat F$. 
Set $\tits H_k=\tau_k^+$, $\tits H_\infty=\tau_\infty^+$ and denote by $\zeta_k$ and $\zeta_\infty$
the respective centers. Because $p_\infty\in h_\infty$ and $H_\infty$ is orthogonal to $\hat F$, the Busemann function $b_{\zeta_\infty,p_\infty}$
is non-negative on $\hat F$. Since the geodesic rays $p_k\zeta_k$ converge to $p_\infty\zeta_\infty$, we find a sequence $t_k\to\infty$ and points $x_k\in p_k\zeta_k$
such that $b_{\zeta_\infty,p_\infty}(x_k)\leq-t_k$. The geodesic ray $p_k\zeta_k$ lies in $H_k$ and therefore $x_k$ lies in $P(h_k)$.
Thus $HB(\zeta_\infty,x_k)\cap P(h_k)$ is a non-empty closed convex set and we have
\[\geo(HB(\zeta_\infty,x_k)\cap P(h_k))=\bar B_{\frac{\pi}{2}}(\zeta_\infty)\cap\geo P(h_k).\]
As $\hat F\subset P(h_k)$, we see $\geo h_\infty\subset \geo(HB(\zeta_\infty,x_k)\cap P(h_k))$.
From Lemma~\ref{lem_flats_in_convex} we infer that there is a $(d-1)$-flat $E_k\subset HB(\zeta_\infty,x_k)\cap P(h_k)$
with $\geo E_k=\geo h_\infty$. Since the flats $h_k$ are pairwise non-parallel and $E_k\subset P(h_k)$,
we see that there is a flat $F_k$ in $P(h_k)$ of dimension at least $d$ which contains $E_k$ and such that $\geo h_k\subset\geo F_k$.
Note  that $b_{\zeta_\infty,p_\infty}$ is bounded above on $E_k$ by $-t_k$ since $E_k\subset HB(\zeta_\infty,x_k)$.
As $b_{\zeta_\infty,p_\infty}$ is non-negative on $\hat F$ we conclude $E_k\cap N_{t_k}(\hat F)=\emptyset$. Hence if $F'_k\subset\hat F$
denotes the closest parallel flat to $F_k$, then $|F_k,F'_k|_H\geq t_k$. In particular, $F_k$ and $F'_k$
bound a flat strip of dimension $\geq d+1$, orthogonal to $\hat F$ and of width at least $t_k$. Since $\hat F$ is periodic and $t_k\to\infty$, we obtain 
a flat half-space $\hat H$ orthogonal to $\hat F$ and of dimension at least $d+1$ as a limit.

To prove the finer statement, simply note that in the argument above the flat $E_k\subset P(h_k)$
spans a flat $F_k\subset P(h_k)$ of dimension $\hat d$ if and only if the convex hull of the spheres $\geo E_k=\geo h_\infty$
and $\geo h_k$ is a round sphere of dimension $\hat d-1$ in $\hat\si$. The rest of the argument stays the same.   
\qed
\medskip

\subsection{Tits boundary of a periodic Morse flat}

In this final section we provide proofs of our main results.

\bdfn
A round sphere $\si\subset\hat\si$ is called {\em singular}, if there exists a round hemisphere $\tau^+\subset\tits X$
with $\D\tau^+=\si$ and  $\tau^+\cap\hat\si=\si$ such that the union $\tau^+\cup\tau^-$ is a round sphere for every 
round hemisphere $\tau^-\subset\hat\si$ with $\D\tau^-=\si$.
We  call a singular sphere $\si\subset\hat\si$ {\em maximal}, if it is not contained in a singular sphere in $\hat\si$ 
of strictly larger dimension.
\edfn

\bcor\label{cor_finite}
Let $X$ be a locally compact  CAT(0) space.
Let $\hat F$ be a periodic $\hat n$-flat in $X$ with $\tits\hat F=\hat \si$.
Then the set of maximal singular spheres in $\hat\si$ is finite.
\ecor

\proof
Proposition~\ref{prop_halfspace} implies that a maximal singular sphere $\si\subset\hat\si$
has to be isolated in the following sense. If $(\si_k)$
is a sequence of singular spheres which converges to a sphere $\si_\infty\subset\si$,
then we must have $\si_k\subset\si$ for almost all $k\in\N$. This yields the claim. 
\qed
\medskip

Recall that  a point in $\hat\si$ is regular, if it has a Tits-neighborhood which is entirely contained in $\hat\si$.
Otherwise we call it {\em singular}.

\blem\label{lem_singptinsingsph}
Every singular point in $\hat\si$ is contained in a singular sphere.
\elem

\proof
If $\xi\in\hat\si$ is a singular point, then there exists a sequence of branch points $\xi_k\in\hat\si$ with $\xi_k\to\xi$.
By Lemma~\ref{lem_perpH}, every such branch point yields a flat half-plane $H_k$ orthogonal to $\hat F$ and with $\xi_k\in\geo H_k$.
Using the cocompact stabilizer of $\hat F$ we can produce a limit flat half-plane $H_\infty$ which is again orthogonal to $\hat F$ and 
with $\xi\in\geo H_\infty$. Set $\tau^+_1:=\geo H_\infty$. If $\D\tau^+_1$ is not singular, then by the Lune Lemma (Lemma~\ref{lem_lune})
we can find a half-circle $\tau^-_1\subset\hat\si$ with $\D\tau^+_1=\D\tau^-_1$ such that $\tau^+_1$ and $\tau^-_1$
span a 2-dimensional spherical lune $\la_2$. By Corollary~\ref{cor_perpH}, there exists a 2-dimensional round hemisphere $\tau^+_2\subset\tits X$
with $\tau^+_2\cap\hat\si=\D\tau^+_2$ and $\tau^-_1\subset\D\tau^+_2$, in particular $\xi\in\D\tau^+_2$. 
If $\D\tau^+_2$ is not singular, we repeat the argument to produce a 3-dimensional round hemisphere 
$\tau^+_3\subset\tits X$ with  $\xi\in\D\tau^+_3$. After at most $(n-1)$ steps we end up with a singular sphere that contains the point $\xi$.
\qed
\medskip

Now we obtain the following refinement of Theorem~\ref{thm_mainF}.

\bthm\label{thm_finitesing}
Let $X$ be a locally compact CAT(0) space.  Suppose that $X$ contains a periodic Morse flat $F$. Then 
the complement in $\tits F$ of the set of regular points can be covered by a finite set of singular spheres of positive codimension.
\ethm

\proof
Immediate from Corollary~\ref{cor_finite} and Lemma~\ref{lem_singptinsingsph}.
\qed

\proof[Proof of Theorem~\ref{thm_mainA}]
Let $F\subset X$ be a periodic $n$-flat.
By Theorem~\ref{thm_mainF}, $\si=\tits F$ contains a  dense subset of spherical points.
Then, by Theorem~\ref{thm_rec_sph}, $\tits X$ is a spherical join or a spherical building.
Since $X$ is geodesically complete and locally compact, the claim follows from Theorem~\ref{thm_leeb} 
and \cite[Proposition~2.3.7]{KleinerLeeb}.
\qed

\proof[Proof of Corollary~\ref{cor_mainA}]
By Theorem~\ref{thm_mainA}, Lemma~\ref{lem_splitflat} and Lemma~\ref{lem_bieber}, we can produce a product decomposition of $X$ as claimed.
Since every such decomposition induces a corresponding join decomposition of $\tits X$.
The uniqueness statement follows from the uniqueness of such join decompositions \cite[Corollary~1.2]{Ly_rigidity}.
\qed

\bibliographystyle{alpha}
\bibliography{rr_I}


\end{document}